\documentclass[reqno,a4paper]{amsart}%
\usepackage{amssymb}
\usepackage[colorlinks=true,hypertex]{hyperref}
\allowdisplaybreaks[4]
\theoremstyle{plain}
\newtheorem{thm}{Theorem}[section]

\newtheorem{conj}{Conjecture}[section]
\theoremstyle{remark}
\newtheorem{rem}{Remark}[section]
\DeclareMathOperator{\td}{d\mspace{-2mu}}
\numberwithin{equation}{section}
\date{Completed on 21 October 2008 at VU's Student Village in Melbourne}
\date{}

\begin{document}

\title[Bounds for the ratio of two gamma functions]
{Bounds for the ratio of two gamma functions---From Gautschi's and Kershaw's inequalities to completely monotonic functions}

\author[F. Qi]{Feng Qi}
\address[F. Qi]{Research Institute of Mathematical Inequality Theory, Henan Polytechnic University, Jiaozuo City, Henan Province, 454010, China}
\email{\href{mailto: F. Qi <qifeng618@gmail.com>}{qifeng618@gmail.com}, \href{mailto: F. Qi <qifeng618@hotmail.com>}{qifeng618@hotmail.com}, \href{mailto: F. Qi <qifeng618@qq.com>}{qifeng618@qq.com}}
\urladdr{\url{http://qifeng618.spaces.live.com}}

\begin{abstract}
In this expository and survey paper, along one of main lines of bounding the ratio of two gamma functions, we look back and analyse some inequalities, the complete monotonicity of several functions involving ratios of two gamma or $q$-gamma functions, the logarithmically complete monotonicity of a function involving the ratio of two gamma functions, some new bounds for the ratio of two gamma functions and divided differences of polygamma functions, and related monotonicity results.
\end{abstract}



\thanks{The author was partially supported by the China Scholarship Council}

\thanks{This paper was typeset using \AmS-\LaTeX}

\maketitle

\tableofcontents

\section{Introduction}

For the sake of proceeding smoothly, we briefly introduce some necessary concepts and notation.

\subsection{The gamma and $q$-gamma functions}

It is well-known that the classical Euler gamma function may be defined by
\begin{equation}\label{egamma}
\Gamma(x)=\int^\infty_0t^{x-1} e^{-t}\td t
\end{equation}
for $x>0$. The logarithmic derivative of $\Gamma(x)$, denoted by $\psi(x)=\frac{\Gamma'(x)}{\Gamma(x)}$, is called the psi or digamma function, and $\psi^{(k)}(x)$ for $k\in \mathbb{N}$ are called the polygamma functions. It is common knowledge that special functions $\Gamma(x)$, $\psi(x)$ and $\psi^{(k)}(x)$ for $k\in\mathbb{N}$ are fundamental and important and have much extensive applications in mathematical sciences.
\par
The $q$-analogues of $\Gamma$ and $\psi$ are defined~\cite[pp.~493\nobreakdash--496]{andrews} for $x>0$ by
\begin{gather}\label{q-gamma-dfn}
\Gamma_q(x)=(1-q)^{1-x}\prod_{i=0}^\infty\frac{1-q^{i+1}}{1-q^{i+x}},\quad 0<q<1,\\
\label{q-gamma-dfn-q>1}
\Gamma_q(x)=(q-1)^{1-x}q^{\binom{x}2}\prod_{i=0}^\infty\frac{1-q^{-(i+1)}}{1-q^{-(i+x)}}, \quad q>1,
\end{gather}
and
\begin{align}\label{q-gamma-1.4}
\psi_q(x)=\frac{\Gamma_q'(x)}{\Gamma_q(x)}&=-\ln(1-q)+\ln q \sum_{k=0}^\infty\frac{q^{k+x}}{1-q^{k+x}}\\
&=-\ln(1-q)-\int_0^\infty\frac{e^{-xt}}{1-e^{-t}}\td\gamma_q(t) \label{q-gamma-1.5}
\end{align}
for $0<q<1$, where $\td\gamma_q(t)$ is a discrete measure with positive masses $-\ln q$ at the positive points $-k\ln q$ for $k\in\mathbb{N}$, more accurately,
\begin{equation}
\gamma_q(t)=
\begin{cases}
-\ln q\sum\limits_{k=1}^\infty\delta(t+k\ln q),&0<q<1,\\ t,&q=1.
\end{cases}
\end{equation}
See~\cite[p.~311]{Ismail-Muldoon-119}.
\par
The $q$-gamma function $\Gamma_q(z)$ has the following basic properties:
\begin{equation}
\lim_{q\to1^+}\Gamma_q(z)=\lim_{q\to1^-}\Gamma_q(z)=\Gamma(z)\quad \text{and}\quad \Gamma_q(x)=q^{\binom{x-1}2}\Gamma_{1/q}(x).
\end{equation}

\subsection{The generalized logarithmic mean}
The generalized logarithmic mean $L_p(a,b)$ of order $p\in\mathbb{R}$ for positive numbers $a$ and $b$ with $a\ne b$ may be defined~\cite[p.~385]{bullenmean} by
\begin{equation}
L_p(a,b)=
\begin{cases}
\left[\dfrac{b^{p+1}-a^{p+1}}{(p+1)(b-a)}\right]^{1/p},&p\ne-1,0;\\[1em]
\dfrac{b-a}{\ln b-\ln a},&p=-1;\\[1em]
\dfrac1e\left(\dfrac{b^b}{a^a}\right)^{1/(b-a)},&p=0.
\end{cases}
\end{equation}
It is well-known that
\begin{gather}
L_{-2}(a,b) =\sqrt{ab}\,=G(a,b),\quad L_{-1}(a,b)=L(a,b),\\
L_0(a,b)=I(a,b)\quad \text{and} \quad L_1(a,b)=\frac{a+b}2=A(a,b)
\end{gather}
are called respectively the geometric mean, the logarithmic mean, the identric or exponential mean, and the arithmetic mean. It is also known~\cite[pp.~386\nobreakdash--387, Theorem~3]{bullenmean} that the generalized logarithmic mean $L_p(a,b)$ of order $p$ is increasing in $p$ for $a\ne b$. Therefore, inequalities
\begin{equation}\label{mean-ineq}
G(a,b)<L(a,b)<I(a,b)<A(a,b)
\end{equation}
are valid for $a>0$ and $b>0$ with $a\ne b$. See also~\cite{abstract-jipam, abstract-rgmia, qi1}. Moreover, the generalized logarithmic mean $L_p(a,b)$ is a special case of $E(r,s;x,y)$ defined by~\eqref{emv-dfn}, that is, $L_p(a,b)=E(1,p+1;a,b)$.

\subsection{Logarithmically completely monotonic functions}

A function $f$ is said to be completely monotonic on an interval $I$ if $f$ has derivatives of all orders on $I$ and
\begin{equation}\label{cmf-dfn-ineq}
(-1)^{n}f^{(n)}(x)\ge0
\end{equation}
for $x \in I$ and $n \ge0$.

\begin{thm}[{\cite[p.~161]{widder}}]\label{p.161-widder}
A necessary and sufficient condition that $f(x)$ should be completely monotonic for $0<x<\infty$ is that
\begin{equation}
f(x)=\int_0^\infty e^{-xt}\td\alpha(t),
\end{equation}
where $\alpha(t)$ is nondecreasing and the integral converges for $0<x<\infty$.
\end{thm}

\begin{thm}[{\cite[p.~83]{bochner}}]\label{p.83-bochner}
If $f(x)$ is completely monotonic on $I$, $g(x)\in I$, and $g'(x)$ is completely monotonic on $(0,\infty)$, then $f(g(x))$ is completely monotonic on $(0,\infty)$.
\end{thm}

A positive function $f(x)$ is said to be logarithmically completely monotonic on an interval $I\subseteq\mathbb{R}$ if it has derivatives of all orders on $I$ and its logarithm $\ln f(x)$ satisfies
$$
(-1)^k[\ln f(x)]^{(k)}\ge0
$$
for $k\in\mathbb{N}$ on $I$.
\par
The notion ``logarithmically completely monotonic function'' was first put forward in~\cite{Atanassov} without an explicit definition. This terminology was explicitly recovered in~\cite{minus-one} whose revised and expanded version was formally published as~\cite{minus-one.tex-rev}.
\par
It has been proved once and again in~\cite{CBerg, clark-ismail-nonlinear, clark-ismail-rgmia, compmon2, absolute-mon.tex, minus-one, minus-one.tex-rev, schur-complete} that a logarithmically completely monotonic function on an interval $I$ must also be completely monotonic on $I$. C. Berg points out in~\cite{CBerg} that these functions are the same as those studied by Horn~\cite{horn} under the name infinitely divisible completely monotonic functions. For more information, please refer to~\cite{CBerg, e-gam-rat-comp-mon, auscm-rgmia} and related references therein.

\subsection{Outline of this paper}

In this expository and survey paper, along one of main lines of bounding the ratio of two gamma functions, we look back and analyse Gautschi's double inequality and Kershaw's second double inequality, the complete monotonicity of several functions involving ratios of two gamma or $q$-gamma functions by Alzer, Bustoz-Ismail, Elezovi\'c-Giordano-Pe\v{c}ari\'c and Ismail-Muldoon, the logarithmically complete monotonicity of a function involving the ratio of two gamma functions, some new bounds for the ratio of two gamma functions and the divided differences of polygamma functions, and related monotonicity results by Batir, Elezovi\'c-Pe\v{c}ari\'c, Qi and others.

\section{Gautschi's and Kershaw's double inequalities}

In this section, we begin with the papers~\cite{gaut, kershaw} to introduce a kind of inequalities for bounding the ratio of two gamma functions.

\subsection{Gautschi's double inequalities}

The first result of the paper~\cite{gaut} was the double inequality
\begin{equation}\label{gaut-3-ineq}
\frac{(x^p+2)^{1/p}-x}2<e^{x^p}\int_x^\infty e^{-t^p}\td t\le c_p\biggl[\biggl(x^p+\frac1{c_p}\biggr)^{1/p}-x\biggr]
\end{equation}
for $x\ge0$ and $p>1$, where
\begin{equation}
c_p=\biggl[\Gamma\biggl(1+\frac1p\biggr)\biggr]^{p/(p-1)}
\end{equation}
or $c_p=1$. By an easy transformation, the inequality~\eqref{gaut-3-ineq} was written in terms of the complementary gamma function
\begin{equation}
\Gamma(a,x)=\int_x^\infty e^{-t}t^{a-1}\td t
\end{equation}
as
\begin{equation}\label{gaut-4-ineq}
\frac{p[(x+2)^{1/p}-x^{1/p}]}2<e^x\Gamma\biggl(\frac1p,x\biggr)\le pc_p\biggl[\biggl(x+\frac1{c_p}\biggr)^{1/p}-x^{1/p}\biggr]
\end{equation}
for $x\ge0$ and $p>1$. In particular, if letting $p\to\infty$, the double inequality
\begin{equation}
\frac12\ln\biggl(1+\frac2x\biggr)\le e^xE_1(x)\le\ln\biggl(1+\frac1x\biggr)
\end{equation}
for the exponential integral $E_1(x)=\Gamma(0,x)$ for $x>0$ was derived from~\eqref{gaut-4-ineq}, in which the bounds exhibit the logarithmic singularity of $E_1(x)$ at $x=0$. As a direct consequence of the inequality~\eqref{gaut-4-ineq} for $p=\frac1s$, $x=0$ and $c_p=1$, the following simple inequality for the gamma function was deduced:
\begin{equation}\label{gaut-none-ineq}
2^{s-1}\le\Gamma(1+s)\le1,\quad 0\le s\le 1.
\end{equation}
\par
The second result of the paper~\cite{gaut} was a sharper and more general inequality
\begin{equation}\label{gaut-6-ineq}
e^{(s-1)\psi(n+1)}\le\frac{\Gamma(n+s)}{\Gamma(n+1)}\le n^{s-1}
\end{equation}
for $0\le s\le1$ and $n\in\mathbb{N}$ than~\eqref{gaut-none-ineq}. It was obtained by proving that the function
\begin{equation}
f(s)=\frac1{1-s}\ln\frac{\Gamma(n+s)}{\Gamma(n+1)}
\end{equation}
is monotonically decreasing for $0\le s<1$ and that
\begin{equation*}
\lim_{s\to1^-}f(s)=-\lim_{s\to1^-}\psi(n+s)=-\psi(n+1).
\end{equation*}

\begin{rem}
For more information on refining the inequality~\eqref{gaut-3-ineq}, please refer to~\cite{incom-gamma-L-N, qi-senlin-mia, Qi-Mei-99-gamma} and related references therein.
\end{rem}

\begin{rem}
The left-hand side inequality in~\eqref{gaut-6-ineq} can be rearranged as
\begin{equation}\label{gaut-ineq-1}
\frac{\Gamma(n+s)}{\Gamma(n+1)}\exp((1-s)\psi(n+1))\ge1
\end{equation}
or
\begin{equation}\label{gaut-ineq-2-exp}
\biggl[\frac{\Gamma(n+s)}{\Gamma(n+1)}\biggr]^{1/(s-1)}e^{-\psi(n+1)}\le1
\end{equation}
for $n\in\mathbb{N}$ and $0\le s\le 1$. Since the limit
\begin{equation}
\lim_{n\to\infty}\biggl\{\biggl[\frac{\Gamma(n+s)}{\Gamma(n+1)}\biggr]^{1/(s-1)} e^{-\psi(n+1)}\biggr\}=1
\end{equation}
can be verified by using Stirling's formula in~\cite[p.~257, 6.1.38]{abram}: For $x>0$, there exists $0<\theta<1$ such that
\begin{equation}\label{Stirling-formula}
\Gamma(x+1)=\sqrt{2\pi}\,x^{x+1/2}\exp\biggl(-x+\frac{\theta}{12x}\biggr),
\end{equation}
it is natural to guess that the function
\begin{equation}\label{gaut-funct}
\biggl[\frac{\Gamma(x+s)}{\Gamma(x+1)}\biggr]^{1/(s-1)}e^{-\psi(x+1)}
\end{equation}
for $0\le s<1$ is possibly increasing with respect to $x$ on $(-s,\infty)$.
\end{rem}

\begin{rem}
For information on the study of the right-hand side inequality in~\eqref{gaut-6-ineq}, please refer to~\cite{bounds-two-gammas.tex, Wendel-Gautschi-type-ineq.tex, Wendel2Elezovic.tex} and a great amount of related references therein.
\end{rem}

\subsection{Kershaw's second double inequality and its proof}\label{kershaw-sec}

In 1983, over twenty years later after the paper~\cite{gaut}, among other things, D.~Kershaw was motivated by the left-hand side inequality~\eqref{gaut-6-ineq} in~\cite{gaut} and presented in~\cite{kershaw} the following double inequality for $0<s<1$ and $x>0$:
\begin{gather}\label{gki2}
\exp\big[(1-s)\psi\big(x+\sqrt{s}\,\big)\big] <\frac{\Gamma(x+1)}{\Gamma(x+s)}
<\exp\biggl[(1-s)\psi\biggl(x+\frac{s+1}2\biggr)\biggr].
\end{gather}
It is called in the literature Kershaw's second double inequality.

\begin{proof}[Kershaw's proof for~\eqref{gki2}]
Define the function $f_\alpha$ by
\begin{equation}\label{kershaw-f-dfn}
  f_\alpha(x)=\frac{\Gamma(x+1)}{\Gamma(x+s)}\exp((s-1)\psi(x+\alpha))
\end{equation}
for $x>0$ and $0<s<1$, where the parameter $\alpha$ is to be determined.
\par
It is not difficult to show, with the aid of Stirling's formula, that
\begin{equation}\label{kershaw-2.3}
  \lim_{x\to\infty}f_\alpha(x)=1.
\end{equation}
\par
Now let
\begin{equation}\label{kershaw-F-dfn}
  F(x)=\frac{f_\alpha(x)}{f_\alpha(x+1)}=\frac{x+s}{x+1}\exp\frac{1-s}{x+\alpha}.
\end{equation}
Then
\begin{equation*}
  \frac{F'(x)}{F(x)}=(1-s)\frac{(\alpha^2-s)+(2\alpha-s-1)x}{(x+1)(x+s)(x+\alpha)^2}.
\end{equation*}
It is easy to show that
\begin{enumerate}
  \item
  if $\alpha=s^{1/2}$, then $F'(x)<0$ for $x>0$;
  \item
  if $\alpha=\frac{s+1}2$, then $F'(x)>0$ for $x>0$.
\end{enumerate}
Consequently if $\alpha=s^{1/2}$ then $F$ strictly decreases, and since $F(x)\to1$ as $x\to\infty$ it follows that $F(x)>1$ for $x>0$. But, from~\eqref{kershaw-2.3}, this implies that $f_\alpha(x)>f_\alpha(x+1)$ for $x>0$, and so $f_\alpha(x)>f_\alpha(x+n)$. Take the limit as $n\to\infty$ to give the result that $f_\alpha(x)>1$, which can be rewritten as the left-hand side inequality in~\eqref{gki2}. The corresponding upper bound can be verified by a similar argument when $\alpha=\frac{s+1}2$, the only difference being that in this case $f_\alpha$ strictly increases to unity.
\end{proof}

\begin{rem}
The idea contained in the above stated proof of~\eqref{gki2} was also utilized by other mathematicians. For detailed information, please refer to related contents and references in~\cite{bounds-two-gammas.tex}.
\end{rem}

\begin{rem}
The inequality~\eqref{gki2} can be rearranged as
\begin{equation}\label{gki2-rew-1}
\frac{\Gamma(x+s)}{\Gamma(x+1)}\exp\big[(1-s)\psi\big(x+\sqrt{s}\,\big)\big] <1 <\frac{\Gamma(x+s)}{\Gamma(x+1)}\exp\biggl[(1-s)\psi\biggl(x+\frac{s+1}2\biggr)\biggr]
\end{equation}
or
\begin{multline}\label{gki2-rew-2}
\biggl[\frac{\Gamma(x+s)}{\Gamma(x+1)}\biggr]^{1/(s-1)}\exp\big[-\psi\big(x+\sqrt{s}\,\big)\big] >1\\
>\biggl[\frac{\Gamma(x+s)}{\Gamma(x+1)}\biggr]^{1/(s-1)}\exp\biggl[-\psi\biggl(x+\frac{s+1}2\biggr)\biggr]. \end{multline}
By Stirling's formula~\eqref{Stirling-formula}, we can prove that
\begin{equation}
\lim_{x\to\infty}\biggl\{\biggl[\frac{\Gamma(x+s)}{\Gamma(x+1)}\biggr]^{1/(s-1)} \exp\big[-\psi\big(x+\sqrt{s}\,\big)\big]\biggr\}=1
\end{equation}
and
\begin{equation}
\lim_{x\to\infty}\biggl\{\biggl[\frac{\Gamma(x+s)}{\Gamma(x+1)}\biggr]^{1/(s-1)} \exp\biggl[-\psi\biggl(x+\frac{s+1}2\biggr)\biggr]\biggr\}=1.
\end{equation}
These clues make us to conjecture that the functions in the very ends of inequalities~\eqref{gki2-rew-1} and~\eqref{gki2-rew-2} are perhaps monotonic with respect to $x$ on $(0,\infty)$.
\end{rem}

\section{Several complete monotonicity results}

The complete monotonicity of the functions in the very ends of inequalities~\eqref{gki2-rew-1} were first demonstrated in~\cite{Bustoz-and-Ismail}, and then several related functions were also proved in~\cite{Alzer1, egp, laj-7.pdf} to be (logarithmically) completely monotonic.

\subsection{Bustoz-Ismail's complete monotonicity results}\label{Bustoz-Ismail-sec}

In 1986, motivated by the double inequality~\eqref{gki2} and other related inequalities, J. Bustoz and M.~E.~H. Ismail revealed in~\cite[Theorem~7 and Theorem~8]{Bustoz-and-Ismail} that
\begin{enumerate}
\item
the function
\begin{equation}\label{bustol-ismail-AM}
\frac{\Gamma(x+s)}{\Gamma(x+1)}\exp\biggl[(1-s)\psi\biggl(x+\frac{s+1}2\biggr)\biggr]
\end{equation}
for $0\le s\le1$ is completely monotonic on $(0,\infty)$; When $0<s<1$, the function~\eqref{bustol-ismail-AM} satisfies $(-1)^nf^{(n)}(x)>0$ for $x>0$;
\item
the function
\begin{equation}\label{bustol-ismail-AMM}
\frac{\Gamma(x+1)}{\Gamma(x+s)}\exp\bigl[(s-1)\psi\bigl(x+s^{1/2}\bigr)\bigr]
\end{equation}
for $0<s<1$ is strictly decreasing on $(0,\infty)$.
\end{enumerate}

\begin{rem}
The proof of the complete monotonicity of the function~\eqref{bustol-ismail-AM} in~\cite[Theorem~7]{Bustoz-and-Ismail} relies on the inequality
\begin{equation}\label{lemma3.1-ism}
(y+a)^{-n}-(y+b)^{-n}>(b-a)n\biggl(y+\frac{a+b}2\biggr)^{-n-1},\quad n>0
\end{equation}
for $y>0$ and $0<a<b$, the series representation
\begin{equation}\label{series-repr}
\psi(x)=-\gamma-\frac1x+\sum_{n=1}^\infty\biggl(\frac1n-\frac1{x+n}\biggr)
\end{equation}
in~\cite[p.~15]{er}, and the above Theorem~\ref{p.83-bochner} applied to $f(x)=e^{-x}$.
\end{rem}

\begin{rem}
The inequality~\eqref{lemma3.1-ism} verified in~\cite[Lemma~3.1]{Bustoz-and-Ismail} can be rewritten as
\begin{equation}\label{lemma3.1-ism-rew}
\biggl[\frac1{-n}\cdot\frac{(y+a)^{-n}-(y+b)^{-n}}{(y+a)-(y+b)}\biggr]^{1/[(-n)-1]}
<\frac{(y+a)+(y+b)}2,\quad n>0
\end{equation}
for $y>0$ and $0<a<b$, which is equivalent to
\begin{equation}\label{E-E(1,2)}
E(1,-n;y+a,y+b)<E(1,2;y+a,y+b),
\end{equation}
where $E(r,s;x,y)$ stands for extended mean values and is defined for two positive numbers $x$ and $y$ and two real numbers $r$ and $s$ by
\begin{equation}
\begin{aligned}\label{emv-dfn}
E(r,s;x,y)&=\biggl(\frac{r}{s}\cdot\frac{y^s-x^s}
{ y^r-x^r}\biggr)^{{1/(s-r)}}, &  rs(r-s)(x-y)&\ne 0; \\
E(r,0;x,y)&=\biggl(\frac{1}{r}\cdot\frac{y^r-x^r}
{\ln y-\ln x}\biggr)^{{1/r}}, &  r(x-y)&\ne 0; \\
E(r,r;x,y)&=\frac1{e^{1/r}}\biggl(\frac{x^{x^r}}{y^{y^r}}\biggr)^{ {1/(x^r-y^r)}},&  r(x-y)&\ne 0; \\
E(0,0;x,y)&=\sqrt{xy}, &  x&\ne y; \\
E(r,s;x,x)&=x, &  x&=y.
\end{aligned}
\end{equation}
Actually, the inequality~\eqref{E-E(1,2)} is an immediate consequence of monotonicity of $E(r,s;x,y)$, see~\cite{ls2}. For more information, please refer to~\cite{bullenmean, ajmaa-mean-chen-qi, emv-log-convex-simple.tex, pqsx, schext-rgmia, schext-rocky, qi1, pams-62, cubo, cubo-rgmia, exp-funct-appl-means-simp.tex, ql, Qi-Luo-1999-Sig, qx1, (b^x-a^x)/x, qx3, zhang-chen-qi-emv} and related references therein.
\end{rem}

\begin{rem}
The proof of the decreasing monotonicity of the function~\eqref{bustol-ismail-AMM} just used the formula~\eqref{series-repr} and and the above Theorem~\ref{p.83-bochner} applied to $f(x)=e^{-x}$.
\end{rem}

\begin{rem}
Indeed, J. Bustoz and M. E. H. Ismail had proved in~\cite[Theorem~7]{Bustoz-and-Ismail} that the function~\eqref{bustol-ismail-AM} is logarithmically completely monotonic on $(0,\infty)$ for $0\le s\le1$. However, because the inequality~\eqref{cmf-dfn-ineq} strictly holds for a completely monotonic function $f$ on $(a,\infty)$ unless $f(x)$ is constant (see~\cite[p.~98]{Dubourdieu}, \cite[p.~82]{e-gam-rat-comp-mon} and~\cite{haerc1}), distinguishing between the cases $0\le s\le1$ and $0<s<1$ is not necessary.
\end{rem}

\subsection{Alzer's and related complete monotonicity results}\label{alzer-comp-sec}

Stimulated by the complete monotonicity obtained in~\cite{Bustoz-and-Ismail}, including those mentioned above, H. Alzer obtained in~\cite[Theorem~1]{Alzer1} that the function
\begin{equation}\label{alzer-func}
\frac{\Gamma(x+s)}{\Gamma(x+1)}\cdot\frac{(x+1)^{x+1/2}}{(x+s)^{x+s-1/2}}
\exp\biggl[s-1+\frac{\psi'(x+1+\alpha)-\psi'(x+s+\alpha)}{12}\biggr]
\end{equation}
for $\alpha>0$ and $s\in(0,1)$ is completely monotonic on $(0,\infty)$ if and only if $\alpha\ge\frac12$, so is the reciprocal of~\eqref{alzer-func} for $\alpha\ge0$ and $s\in(0,1)$ if and only if $\alpha=0$.
\par
As consequences of the monotonicity of the function~\eqref{alzer-func}, the following inequalities are deduced in~\cite[Corollary~1 and Corollary~2]{Alzer1}:
\begin{enumerate}
\item
The inequalities
\begin{equation}\label{alzer-fun-ineq}
\begin{gathered}
\exp\biggl[s-1+\frac{\psi'(x+1+\beta) -\psi'(x+s+\beta)}{12}\biggr] \le\frac{(x+s)^{x+s-1/2}}{(x+1)^{x+1/2}}\cdot\frac{\Gamma(x+1)}{\Gamma(x+s)}\\
\le \exp\biggl[s-1+\frac{\psi'(x+1+\alpha) -\psi'(x+s+\alpha)}{12}\biggr],\quad \alpha>\beta\ge0
\end{gathered}
\end{equation}
are valid for all $s\in(0,1)$ and $x\in(0,\infty)$ if and only if $\beta=0$ and $\alpha\ge\frac12$.
\item
If
\begin{equation}
a_n=\frac32\biggl\{1+\ln\biggl[\frac{2[\Gamma((n+1)/2)]^2} {[\Gamma(n/2)]^2}\cdot{n^{n-1}}{(n+1)^n}\biggr]\biggr\},
\end{equation}
then
\begin{equation}\label{sum-alzer-ineq}
a_n<(-1)^{n+1}\Biggl[\frac{\pi^2}{12}-\sum_{k=1}^n(-1)^{k+1}\frac1{k^2}\Biggr]<a_{n+1},\quad n\in\mathbb{N}.
\end{equation}
\end{enumerate}

\begin{rem}
The inequality~\eqref{sum-alzer-ineq} follows from the formula
\begin{equation}
\frac14\biggl[\psi'\biggl(\frac{n}2+1\biggr)-\psi'\biggl(\frac{n+1}2\biggr)\biggr] =\sum_{k=1}^\infty\frac{(-1)^k}{(n+k)^2} =(-1)^{n}\Biggl[\sum_{k=1}^n\frac{(-1)^{k+1}}{k^2}-\frac{\pi^2}2\Biggr]
\end{equation}
and the inequality~\eqref{alzer-fun-ineq} applied to $s=\frac12$, $\alpha=\frac12$ and $\beta=0$.
\end{rem}

\begin{rem}
The proof of the complete monotonicity of the function~\eqref{alzer-func} in~\cite{Alzer1} is based on Theorem~\ref{p.83-bochner} applied to $f(x)=e^{-x}$, the formulas
\begin{equation}
\frac1x=\int_0^\infty e^{-xt}\td t,\quad \ln\frac{y}x=\int_0^\infty\frac{e^{-xt}-e^{-yt}}t\td t
\end{equation}
and
\begin{equation}
\psi(x)=-\gamma+\int_0^\infty\frac{e^{-t}-e^{-xt}}{1-e^{-t}}\td t
\end{equation}
for $x,y>0$, and discussing the positivity of the functions
\begin{equation}
\frac{12-t^2e^{-\alpha t}}{12(1-e^{-t})}-\frac12-\frac1t\quad\text{and}\quad
\frac12+\frac1t-\frac{12-t^2}{12(1-e^{-t})}
\end{equation}
for $x\in(0,\infty)$ and $\alpha\ge\frac12$. Therefore, H. Alzer essentially gave in~\cite[Theorem~1]{Alzer1} necessary and sufficient conditions for the function~\eqref{alzer-func} to be logarithmically completely monotonic on $(0,\infty)$.
\end{rem}

\begin{rem}
In~\cite[Theorem~3]{laj-7.pdf}, a slight extension of~\cite[Theorem~1]{Alzer1} was presented: The function
\begin{equation}\label{li-ext-fun}
\frac{\Gamma(x+s)}{\Gamma(x+t)}\cdot\frac{(x+t)^{x+t-1/2}}{(x+s)^{x+s-1/2}}
\exp\biggl[s-t+\frac{\psi'(x+t+\alpha) -\psi'(x+s+\alpha)}{12}\biggr]
\end{equation}
for $0<s<t$ and $x\in(0,\infty)$ is logarithmically completely monotonic if and only if $\alpha\ge\frac12$, so is the reciprocal of~\eqref{li-ext-fun} if and only if $\alpha=0$.
\par
The decreasing monotonicity of~\eqref{li-ext-fun} and its reciprocal imply that the double inequality
\begin{multline}\label{li-ext-fun-ineq}
\exp\biggl[t-s+\frac{\psi'(x+s+\beta) -\psi'(x+t+\beta)}{12}\biggr] \le\frac{(x+t)^{x+t-1/2}}{(x+s)^{x+s-1/2}}\cdot\frac{\Gamma(x+s)}{\Gamma(x+t)}\\
\le \exp\biggl[t-s+\frac{\psi'(x+s+\alpha) -\psi'(x+t+\alpha)}{12}\biggr]
\end{multline}
for $\alpha>\beta\ge0$ are valid for $0<s<t$ and $x\in(0,\infty)$ if and only if $\beta=0$ and $\alpha\ge\frac12$.
\par
It is obvious that the inequality~\eqref{li-ext-fun-ineq} is a slight extension of the double inequality~\eqref{alzer-fun-ineq} obtained in~\cite[Corollary~2]{Alzer1}.
\end{rem}

\begin{rem}
In~\cite[Theorem~3.4]{Ismail-Muldoon-119}, the following complete monotonicity were established: Let $0<q<1$ and
\begin{equation}\label{g{alpha,q}(x)}
g_{\alpha,q}(x)=(1-q)^x(1-q^x)^{1/2}\Gamma_q(x) \exp\biggl[\frac{F(q^x)}{\ln q}-\frac{\psi_q'(x+\alpha)}{12}\biggr],
\end{equation}
where
\begin{equation}
F(x)=\sum_{n=1}^\infty\frac{x^n}{n^2}=-\int_0^x\frac{\ln(1-t)}{t}\td t.
\end{equation}
Then $[\ln g_{\alpha,q}(x)]'$ is completely monotonic on $(0,\infty)$ for $\alpha\ge\frac12$, $-[\ln g_{\alpha,q}(x)]'$ is completely monotonic on $(0,\infty)$ for $\alpha\le0$, and neither is completely monotonic on $(0,\infty)$ for $0<\alpha<\frac12$.
\par
As a consequence of~\cite[Theorem~3.4]{Ismail-Muldoon-119}, the following result was deduced in~\cite[Corollary~3.5]{Ismail-Muldoon-119}: Let $0<q<1$, $0<s<1$ and
\begin{equation}
\begin{split}\label{f-alpha-q}
f_{\alpha,q}(x)&=\frac{g_\alpha(x+s)}{g_\alpha(x+1)}\\
&=\frac{(1-q)^{s-1}(1-q^{x+s})^{1/2}\Gamma_q(x+s)}{(1-q^{x+1})^{1/2}\Gamma_q(x+1)}\\ &\quad\times\exp\biggl[\frac{F(q^{x+s})-F(q^{x+1})}{\ln q}+\frac{\psi_q'(x+1+\alpha)-\psi_q'(x+s+\alpha)}{12}\biggr].
\end{split}
\end{equation}
Then $[\ln f_{\alpha,q}(x)]'$ is completely monotonic on $(0,\infty)$ for $\alpha\ge\frac12$, $-[\ln f_{\alpha,q}(x)]'$ is complete monotonic on $(0,\infty)$ for $\alpha\le0$, and neither is completely monotonic on $(0,\infty)$ for $0<\alpha<\frac12$.
\par
Taking the limit $q\to1^-$ in~\eqref{f-alpha-q} yields~\cite[Corollary~3.6]{Ismail-Muldoon-119}, a recovery of~\cite[Theorem~1]{Alzer1} mentioned above.
\end{rem}

\begin{rem}
It is clear that~\cite[Theorem~3]{laj-7.pdf} can be derived by taking the limit
\begin{equation}
\lim_{q\to1^-}\frac{g_\alpha(x+s)}{g_\alpha(x+t)}
\end{equation}
for $0<s<t$, where $g_\alpha(x)$ is defined by~\eqref{g{alpha,q}(x)}.
\end{rem}

\subsection{Ismail-Muldoon's complete monotonicity results}

Inspired by inequalities~\eqref{gaut-6-ineq} and~\eqref{gki2}, Ismail and Muldoon proved in~\cite[Theorem~3.2]{Ismail-Muldoon-119} the following conclusions: For $0<a<b$ and $0<q<1$, let
\begin{equation}\label{Gamma-q(x+a)}
h(x)=\ln\biggl\{\frac{\Gamma_q(x+a)}{\Gamma_q(x+b)}\exp[(b-a)\psi_q(x+c)]\biggr\}.
\end{equation}
If $c\ge\frac{a+b}2$, then $-h'(x)$ is completely monotonic on $(-a,\infty)$; If $c\le a$, then $h'(x)$ is completely monotonic on $(-c,\infty)$; Neither $h'(x)$ or $-h'(x)$ is completely monotonic for $a<c<\frac{a+b}2$. Consequently, the following inequality was deduced in~\cite[Theorem~3.3]{Ismail-Muldoon-119}: If $0<q<1$, the inequality
\begin{equation}\label{q(x+1)}
\frac{\Gamma_q(x+1)}{\Gamma_q(x+s)} <\exp\biggl[(1-s)\psi_q\biggl(x+\frac{s+1}2\biggr)\biggr],\quad 0<s<1
\end{equation}
holds for $x>-s$.
\par
Influenced by~\eqref{q(x+1)}, H. Alzer posed in the final of the paper~\cite[p.~13]{Alzer-Math-Nachr-2001} the following open problem: For real numbers $0<q\ne1$ and $s\in(0,1)$, determine the best possible values $a(q,s)$ and $b(q,s)$ such that the inequalities
\begin{equation}
\exp[(1-s)\psi_q(x+a(q,s))]<\frac{\Gamma_q(x+1)}{\Gamma_q(x+s)} <\exp[(1-s)\psi_q(x+b(q,s))]
\end{equation}
hold for all $x>0$.

\begin{rem}
Since the paper~\cite{Ismail-Muldoon-119} was published in a conference proceedings, it is not easy to acquire it, so the completely monotonic properties of the function $h(x)$, obtained in~\cite[Theorem~3.2]{Ismail-Muldoon-119}, were neglected in most circumstances.
\end{rem}

\subsection{Elezovi\'c-Giordano-Pe\v{c}ari\'c's inequality and monotonicity results}

Inspired by the double inequality~\eqref{gki2}, the following problem was posed in~\cite[p.~247]{egp}: What are the best constants $\alpha$ and $\beta$ such that the double inequality
\begin{equation}
\psi(x+\alpha)\le\frac1{t-s}\int_s^t\psi(u)\td u\le\psi(x+\beta)
\end{equation}
holds for $x>-\min\{s,t,\alpha,\beta\}$?
\par
An answer to the above problem was procured in~\cite[Theorem~4]{egp}:  The double inequality
\begin{equation}\label{second-egp-thm4}
\psi\biggl(x+\psi^{-1}\biggl(\frac1{t-s}\int_s^t\psi(u)\td u\biggr)\biggr)
<\frac1{t-s}\int_s^t\psi(x+u)\td u<\psi\biggl(x+\frac{s+t}2\biggr)
\end{equation}
is valid for every $x\ge0$ and positive numbers $s$ and $t$.
\par
Moreover, the function
\begin{equation}\label{gamma-arithmetic-funct}
\psi\biggl(x+\frac{s+t}2\biggr)-\frac1{t-s}\ln\frac{\Gamma(x+t)}{\Gamma(x+s)}
\end{equation}
for $s,t>0$ and $r=\min\{s,t\}$ was proved in~\cite[Theorem~5]{egp} to be completely monotonic on $(-r,\infty)$.

\begin{rem}
It is clear that~\cite[Theorem~5]{egp} stated above extends or generalizes the complete monotonicity of the function~\eqref{bustol-ismail-AM}.
\end{rem}

\begin{rem}
By the way, the complete monotonicity in~\cite[Theorem~5]{egp} was extended and iterated in~\cite[Proposition~5]{notes-best.tex-mia} and~\cite[Proposition~5]{notes-best.tex-rgmia} as follows: The function
\begin{equation}\label{cmf-lcmf}
\bigg[\dfrac{\Gamma(x+t)}{\Gamma(x+s)}\bigg]^{1/(s-t)} \exp\biggl[\psi\biggl(x+\frac{s+t}2\biggr)\biggr]
\end{equation}
is logarithmically completely monotonic with respect to $x$ on $(-\alpha,\infty)$, where $s$ and $t$ are real numbers and $\alpha=\min\{s,t\}$.
\end{rem}

\begin{rem}
Along the same line as proving the inequality~\eqref{second-egp-thm4} in~\cite{egp}, the inequality~\eqref{second-egp-thm4} was generalized in \cite[Theorem~2]{Chen-Ai-Jun-rgmia-07} as
\begin{multline}\label{chen-ai-jun-rgmia-07-ineq}
(-1)^n\psi^{(n)}\biggl(x+\bigl(\psi^{(n)}\bigr)^{-1}\biggl(\frac1{t-s}\int_s^t\psi^{(n)}(u)\td u\biggr)\biggr)<\\ \frac{(-1)^n\bigl[\psi^{(n-1)}(x+t)-\psi^{(n-1)}(x+s)\bigr]}{t-s} <(-1)^n\psi^{(n)}\biggl(x+\frac{s+t}2\biggr)
\end{multline}
for $x>0$, $n\ge0$, and $s,t>0$, where $\bigl(\psi^{(n)}\bigr)^{-1}$ denotes the inverse function of $\psi^{(n)}$.
\end{rem}

\begin{rem}
Since the inverse functions of the psi and polygamma functions are involved, it is much difficult to calculate the lower bounds in~\eqref{second-egp-thm4} and~\eqref{chen-ai-jun-rgmia-07-ineq}.
\end{rem}

\begin{rem}
In~\cite{kershaw-anal.appl}, by the method used in~\cite{kershaw}, it was proved that the double inequality
\begin{equation}\label{kershaw-singapore-ineq}
\psi\bigl(x+\sqrt{st}\,\bigr)<\frac{\ln\Gamma(x+t)-\ln\Gamma(x+s)}{t-s}<\psi\biggl(x+\frac{s+t}2\biggr)
\end{equation}
holds for $s,t>0$.
It s clear that the upper bound in~\eqref{kershaw-singapore-ineq} is a recovery of~\eqref{second-egp-thm4} and an immediate consequence of the complete monotonicity of the function~\eqref{gamma-arithmetic-funct}.
\end{rem}

\section{Two logarithmically complete monotonicity results}

Suggested by the double inequality~\eqref{gki2}, it is natural to put forward the following problem: What are the best constants $\delta_1(s,t)$ and $\delta_2(s,t)$ such that
\begin{equation}\label{gki2-gen}
\exp[\psi(x+\delta_1(s,t))] \le\biggl[\frac{\Gamma(x+t)}{\Gamma(x+s)}\biggr]^{1/(t-s)} \le\exp[\psi(x+\delta_2(s,t))]
\end{equation}
is valid for $x>-\min\{s,t,\delta_1(s,t),\delta_2(s,t)\}$? where $s$ and $t$ are real numbers.
\par
It is clear that the inequality~\eqref{gki2-gen} can also be rewritten as
\begin{equation}\label{gki2-gen-rew-1}
\biggl[\frac{\Gamma(x+t)}{\Gamma(x+s)}\biggr]^{1/(s-t)}\exp[\psi(x+\delta_1)] \le1 \le\biggl[\frac{\Gamma(x+t)}{\Gamma(x+s)}\biggr]^{1/(s-t)} \exp[\psi(x+\delta_2)]
\end{equation}
which suggests some monotonic properties of the function
\begin{equation}\label{gamma-delta-ratio}
\biggl[\frac{\Gamma(x+t)}{\Gamma(x+s)}\biggr]^{1/(t-s)}\exp[-\psi(x+\delta(s,t))],
\end{equation}
since the limit of the function~\eqref{gamma-delta-ratio} as $x\to\infty$ is $1$ by using~\eqref{Stirling-formula}.
\par
This problem was considered in~\cite{ratio-gamma-polynomial.tex-jcam, ratio-gamma-polynomial.tex-rgmia, gamma-batir.tex-jcam, gamma-batir.tex-rgmia} along two different approaches and the following results of different forms were established.

\begin{thm}[{\cite[Theorem~1]{ratio-gamma-polynomial.tex-jcam} and~\cite[Theorem~1]{ratio-gamma-polynomial.tex-rgmia}}]\label{gamma-ratio-multply}
Let $a$, $b$, $c$ be real numbers and $\rho=\min\{a,b,c\}$. Define
\begin{equation}\label{gamma-multply}
F_{a,b;c}(x)=
\begin{cases}
\biggl[\dfrac{\Gamma(x+b)}{\Gamma(x+a)}\biggr]^{1/(a-b)}\exp[\psi(x+c)], &a\ne
b\\
\exp[\psi(x+c)-\psi(x+a)],&a=b\ne c
\end{cases}
\end{equation}
for $x\in(-\rho,\infty)$. Furthermore, let $\theta(t)$ be an implicit function
defined by equation
\begin{equation}\label{implicit}
e^t-t =e^{\theta(t)}-\theta(t)
\end{equation}
on $(-\infty,\infty)$. Then $\theta(t)$ is decreasing and $t\theta(t)<0$ for
$\theta(t)\ne t$, and
\begin{enumerate}
\item
$F_{a,b;c}(x)$ is logarithmically completely monotonic on $(-\rho,\infty)$ if
\begin{equation}\label{d1}
\begin{split}
(a,b;c)&\in \{c\ge a,c\ge b\}\cup\{c\ge a,0\ge c-b\ge\theta(c-a)\}\\*
&\quad\cup\{c\le a,c-b\ge\theta(c-a)\}\setminus\{a=b=c\};
\end{split}
\end{equation}
\item
$[F_{a,b;c}(x)]^{-1}$ is logarithmically completely monotonic on $(-\rho,\infty)$ if
\begin{equation}\label{d2}
\begin{split}
(a,b;c)&\in \{c\le a,c\le b\}\cup\{c\ge a,c-b\le\theta(c-a)\}\\*
&\quad\cup\{c\le a,0\le c-b\le\theta(c-a)\}\setminus\{a=b=c\}.
\end{split}
\end{equation}
\end{enumerate}
\end{thm}

\begin{thm}[{\cite[Theorem~1]{gamma-batir.tex-jcam} and~\cite[Theorem~1]{gamma-batir.tex-rgmia}}]\label{nu-log-mon}
For real numbers $s$ and $t$ with $s\ne t$ and $\theta(s,t)$ a constant depending on $s$ and $t$, define
\begin{equation}\label{nudef}
\nu_{s,t}(x)=\frac1{\exp\big[\psi\bigl(x+\theta(s,t)\bigr)\big]}
\biggl[\frac{\Gamma(x+t)}{\Gamma(x+s)}\biggr]^{1/(t-s)}.
\end{equation}
\begin{enumerate}
\item
The function $\nu_{s,t}(x)$ is logarithmically completely monotonic on the interval $(-\theta(s,t),\infty)$ if and only if $\theta(s,t)\le\min\{s,t\}$;
\item
The function $[\nu_{s,t}(x)]^{-1}$ is logarithmically completely monotonic on the interval $(-\min\{s,t\},\infty)$ if and only if $\theta(s,t)\ge\frac{s+t}2$.
\end{enumerate}
\end{thm}

\begin{rem}
In~\cite{ratio-gamma-polynomial.tex-jcam, ratio-gamma-polynomial.tex-rgmia}, it was deduced by standard argument that
\begin{gather*}
(-1)^i[\ln F_{a,b;c}(x)]^{(i)}
=\int_0^\infty\biggl[\frac{e^{(c-a)u}-e^{(c-b)u}}{u(b-a)}-1\biggr]
\frac{u^ie^{-(x+c)u}}{1-e^{-u}}\td u\\
=\int_0^\infty\biggl[\frac{[e^{(c-a)u}-(c-a)u]-[e^{(c-b)u}-(c-b)u]} {[(c-a)-(c-b)]u}\biggr] \frac{u^ie^{-(x+c)u}}{1-e^{-u}}\td u
\end{gather*}
for $i\in\mathbb{N}$ and $a\ne b$. Therefore, the sufficient conditions in~\cite[Theorem~1]{ratio-gamma-polynomial.tex-jcam} and~\cite[Theorem~1]{ratio-gamma-polynomial.tex-rgmia} are stated in terms of the implicit function $\theta(t)$ defined by~\eqref{implicit}.
\end{rem}

\begin{rem}
In~\cite{gamma-batir.tex-jcam, gamma-batir.tex-rgmia}, the logarithmic derivative of $\nu_{s,t}(x)$ was rearranged as
\begin{equation}
\ln\nu_{s,t}(x)=\int_0^\infty\frac{e^{-[x+\theta(s,t)]u}}{1-e^{-u}} \Bigl\{1-e^{u[\theta(s,t)+\ln p_{s,t}(u)]}\Bigr\}\td u,
\end{equation}
where
\begin{equation}
p_{s,t}(u)=\biggl(\frac1{t-s}\int_s^t e^{-uv}\td v\biggr)^{1/u}.
\end{equation}
Since the function $p_{s,t}(u)$ is increasing on $[0,\infty)$ with
\begin{equation}
\lim_{u\to0}p_{s,t}(u)=e^{-(s+t)/2}\quad \text{and}\quad \lim_{u\to\infty}p_{s,t}(u)=e^{-\min\{s,t\}},
\end{equation}
the necessary and sufficient conditions in~\cite[Theorem~1]{gamma-batir.tex-jcam} and~\cite[Theorem~1]{gamma-batir.tex-rgmia} may be derived immediately by considering Theorem~\ref{p.161-widder}.
\par
However, the necessary conditions in~\cite[Theorem~1]{gamma-batir.tex-jcam} and~\cite[Theorem~1]{gamma-batir.tex-rgmia} were proved by establishing the following inequalities involving the polygamma functions and their inverse functions in~\cite[Proposition~1]{gamma-batir.tex-jcam} and~\cite[Proposition~1]{gamma-batir.tex-rgmia}:
\begin{enumerate}
\item
If $m>n\ge0$ are two integers, then
\begin{equation}\label{m>n>0}
\left(\psi^{(m)}\right)^{-1}\left(\frac1{t-s} \int_s^t\psi^{(m)}(v)\td
v\right) \le\left(\psi^{(n)}\right)^{-1}\left(\frac1{t-s}
\int_s^t\psi^{(n)}(v)\td v\right),
\end{equation}
where $\left(\psi^{(k)}\right)^{-1}$ stands for the inverse function of $\psi^{(k)}$ for $k\ge0$;
\item
The inequality
\begin{equation}\label{log-mean-ineq}
\psi^{(i)}(L(s,t))\le\frac1{t-s}\int^t_s\psi^{(i)}(u)\td u
\end{equation}
is valid for $i$ being positive odd number or zero and reversed for $i$ being positive even number;
\item
The function
\begin{equation}\label{increas-conc}
\left(\psi^{(\ell)}\right)^{-1}\left(\frac1{t-s} \int_s^t\psi^{(\ell)}(x+v)\td
v\right)-x
\end{equation}
for $\ell\ge0$ is increasing and concave in $x>-\min\{s,t\}$ and has a sharp
upper bound $\frac{s+t}{2}$.
\end{enumerate}
\par
Note that if taking $m=1$, $n=0$, $i=0$ and $\ell=0$ in~\eqref{m>n>0}, \eqref{log-mean-ineq} and~\eqref{increas-conc}, then~\cite[Lemma~1]{f-mean} and~\cite[Theorem~6]{f-mean} may be derived straightforwardly.
\end{rem}

\section{New bounds and monotonicity results}

\subsection{Elezovi\'c-Pe\v{c}ari\'c's lower bound}
The inequality~\eqref{log-mean-ineq} for $i=0$, that is, \cite[Lemma~1]{f-mean}, may be rewritten as
\begin{equation}
\frac{\ln\Gamma(t)-\ln\Gamma(s)}{t-s}\ge\psi(L(s,t))
\end{equation}
or
\begin{equation}\label{Elezovic-Pecaric-ineq-lower}
\biggl[\frac{\Gamma(t)}{\Gamma(s)}\biggr]^{1/(t-s)}\ge e^{\psi(L(s,t))}
\end{equation}
for positive numbers $s$ and $t$.

\begin{rem}
From the left-hand side inequality in~\eqref{mean-ineq}, it is easy to see that the inequality~\eqref{Elezovic-Pecaric-ineq-lower} refines the traditionally lower bound $e^{\psi(G(s,t))}$.
\end{rem}

\begin{rem}
In~\cite[Theorem~2.4]{gamma-fun-ineq-batir}, the following incorrect double inequality
was obtained:
\begin{equation}\label{kershaw-batir}
e^{(x-y)\psi(L(x+1,y+1)-1)} \le\frac{\Gamma(x)}{\Gamma(y)} \le e^{(x-y)\psi(A(x,y))},
\end{equation}
where $x$ and $y$ are positive real numbers.
\end{rem}

\subsection{Allasia-Giordano-Pe\v{c}ari\'c's inequalities}

In Section~4 of \cite{Allasia-Gior-Pecaric-MIA-02}, as straightforward consequences of Hadamard type inequalities obtained in~\cite{agpit}, the following double inequalities for bounding $\ln\frac{\Gamma(y)}{\Gamma(x)}$ were listed: For $y>x>0$, $n\in\mathbb{N}$ and $h=\frac{y-x}n$, we have
\begin{gather*}
\frac{h}2[\psi(x)+\psi(y)]+h\sum_{k=1}^{n-1}\psi(x+kh)<\ln\frac{\Gamma(y)}{\Gamma(x)} <h\sum_{k=0}^{n-1}\psi\biggl(x+\biggl(k+\frac12\biggr)h\biggr),\\
0<h\sum_{k=0}^{n-1}\psi\biggl(x+\biggl(k+\frac12\biggr)h\biggr) -\ln\frac{\Gamma(y)}{\Gamma(x)}
<\ln\frac{\Gamma(y)}{\Gamma(x)} -\frac{h}2[\psi(x)+\psi(y)]-h\sum_{k=1}^{n-1}\psi(x+kh),\\
\frac{h}2[\psi(x)+\psi(y)]+h\sum_{k=1}^{n-1}\psi(x+kh) -\sum_{i=1}^{m-1}\frac{B_{2i}h^{2i}}{(2i)!}\bigl[\psi^{(2i-1)}(y)-\psi^{(2i-1)}(x)\bigr]< \\
\ln\frac{\Gamma(y)}{\Gamma(x)}
<h\sum_{k=0}^{n-1}\psi\biggl(x+\biggl(k+\frac12\biggr)h\biggr) -\sum_{i=1}^{m-1}\frac{B_{2i}h^{2i}}{(2i)!}\bigl[\psi^{(2i-1)}(y)-\psi^{(2i-1)}(x)\bigr],\\
0<h\sum_{k=0}^{n-1}\psi\biggl(x+\biggl(k+\frac12\biggr)h\biggr)  -\sum_{i=1}^{m-1}\frac{B_{2i}(1/2)h^{2i}}{(2i)!}\bigl[\psi^{(2i-1)}(y)-\psi^{(2i-1)}(x)\bigr]
-\ln\frac{\Gamma(y)}{\Gamma(x)}\\
<\ln\frac{\Gamma(y)}{\Gamma(x)}-\frac{h}2[\psi(x)+\psi(y)]-h\sum_{k=1}^{n-1}\psi(x+kh)
+\sum_{i=1}^{m-1}\frac{B_{2i}h^{2i}}{(2i)!}\bigl[\psi^{(2i-1)}(y)-\psi^{(2i-1)}(x)\bigr],\\
h\sum_{k=0}^{n-1}\psi\biggl(x+\biggl(k+\frac12\biggr)h\biggr) -\sum_{i=1}^{m-2}\frac{B_{2i}(1/2)h^{2i}}{(2i)!}\bigl[\psi^{(2i-1)}(y)-\psi^{(2i-1)}(x)\bigr]
<\ln\frac{\Gamma(y)}{\Gamma(x)}\\
<h\sum_{k=0}^{n-1}\psi\biggl(x+\biggl(k+\frac12\biggr)h\biggr) -\sum_{i=1}^{m-1}\frac{B_{2i}(1/2)h^{2i}}{(2i)!}\bigl[\psi^{(2i-1)}(y)-\psi^{(2i-1)}(x)\bigr],\\
\frac{h}2[\psi(x)+\psi(y)]+h\sum_{k=1}^{n-1}\psi(x+kh) -\sum_{i=1}^{m-1}\frac{B_{2i}h^{2i}}{(2i)!}\bigl[\psi^{(2i-1)}(y)-\psi^{(2i-1)}(x)\bigr]
<\ln\frac{\Gamma(y)}{\Gamma(x)}\\
<\frac{h}2[\psi(x)+\psi(y)]+h\sum_{k=1}^{n-1}\psi(x+kh) -\sum_{i=1}^{m-2}\frac{B_{2i}h^{2i}}{(2i)!}\bigl[\psi^{(2i-1)}(y)-\psi^{(2i-1)}(x)\bigr],
\end{gather*}
where $m$ is an odd and positive integer,
\begin{equation}
B_k\biggl(\frac12\biggr)=\biggl(\frac1{2^{k-1}}-1\biggr)B_k,\quad k\ge0
\end{equation}
and $B_i$ for $i\ge0$ are Bernoulli numbers defined by
\begin{equation}
\frac{t}{e^t-1}=\sum_{i=0}^\infty B_i\frac{t^i}{i!} =1-\frac{x}2+\sum_{j=1}^\infty B_{2j}\frac{x^{2j}}{(2j)!}, \quad\vert x\vert<2\pi.
\end{equation}
If replacing $m$ by an even and positive integer, then the last four double inequalities are reversed.

\subsection{Batir's double inequality for polygamma functions}

It is clear that the double inequality~\eqref{gki2} can be rearranged as
\begin{equation}\label{kershaw-rearr}
\psi\bigl(x+\sqrt{s}\,\bigr)<\frac{\ln\Gamma(x+1)-\ln\Gamma(x+s)}{1-s} <\psi\biggl(x+\frac{s+1}2\biggr)
\end{equation}
for $0<s<1$ and $x>1$. The middle term in~\eqref{kershaw-rearr} can be regarded as a divided difference of the function $\ln\Gamma(t)$ on $(x+s,x+1)$. Stimulated by this, N.~Batir extended and generalized in~\cite[Theorem~2.7]{batir-jmaa-06-05-065} the double inequality~\eqref{kershaw-rearr} as
\begin{equation}\label{batir-psi-ineq}
-\bigl\vert \psi^{(n+1)}(L_{-(n+2)}(x,y))\bigr\vert <\frac{\bigl\vert \psi^{(n)}(x)\bigr\vert  -\bigl\vert \psi^{(n)}(y)\bigr\vert }{x-y}<-\bigl\vert \psi^{(n+1)}(A(x,y))\bigr\vert
\end{equation}
where $x,y$ are positive numbers and $n\in\mathbb{N}$.

\subsection{Chen's double inequality in terms of polygamma functions}
In~\cite[Theorem~2]{chen-mean-GK-ineq}, by virtue of the composite Simpson rule
\begin{equation}
\int_a^bf(t)\td t=\frac{b-a}6\biggl[f(a)+4f\biggl(\frac{a+b}2\biggr)+f(b)\biggr] -\frac{(b-a)^5}{2880}f^{(4)}(\xi),\quad \xi\in(a,b)
\end{equation}
in~\cite{Hammerlin-Hoffmann-bbok-91} and the formula
\begin{equation}
\frac1{y-x}\int_x^yf(t)\td t=\sum_{k=0}^\infty\frac1{(2k+1)!} \biggl(\frac{y-x}2\biggr)^{2k}f^{(2k)}\biggl(\frac{x+y}2\biggr)
\end{equation}
in~\cite{Neuman-Sandor-05-Aus}, the following double inequalities and series representations were trivially shown: For $n\in\mathbb{N}$ and positive numbers $x$ and $y$ with $x\ne y$,
\begin{gather*}
\frac13A(\psi(x),\psi(y))+\frac23\psi(A(x,y))-\frac{(y-x)^4}{2880}\psi^{(4)}(\max\{x,y\}) <\frac{\ln\Gamma(y)-\ln\Gamma(x)}{y-x} \\
<\frac13A(\psi(x),\psi(y))+\frac23\psi(A(x,y))-\frac{(y-x)^4}{2880}\psi^{(4)}(\min\{x,y\}),\\
(-1)^{n-1}\biggl[\frac{A\bigl(\psi^{(n)}(x),\psi^{(n)}(y)\bigr)}3 +\frac{2\psi^{(n)}(A(x,y))}3-\frac{(y-x)^4\psi^{(n+4)}(\min\{x,y\})}{2880}\biggr]\\
<\frac{(-1)^{n-1}\bigl[\psi^{(n-1)}(y)-\psi^{(n-1)}(x)\bigr]}{y-x} \\
<(-1)^{n-1}\biggl[\frac{A\bigl(\psi^{(n)}(x),\psi^{(n)}(y)\bigr)}3 +\frac{2\psi^{(n)}(A(x,y))}3-\frac{(y-x)^4\psi^{(n+4)}(\max\{x,y\})}{2880}\biggr],\\
\frac{\ln\Gamma(y)-\ln\Gamma(x)}{y-x}=\sum_{k=0}^\infty\frac1{(2k+1)!} \biggl(\frac{y-x}2\biggr)^{2k}\psi^{(2k)}\biggl(\frac{x+y}2\biggr),\\
\frac{\psi^{(n-1)}(y)-\psi^{(n-1)}(x)}{y-x}=\sum_{k=0}^\infty\frac1{(2k+1)!} \biggl(\frac{y-x}2\biggr)^{2k}\psi^{(2k+n)}\biggl(\frac{x+y}2\biggr).
\end{gather*}

\subsection{Recent monotonicity results by Qi and his coauthors}

Motivated by the left-hand side inequality in~\eqref{kershaw-batir}, although it is not correct, several refinements and generalizations about inequalities~\eqref{Elezovic-Pecaric-ineq-lower} and~\eqref{batir-psi-ineq} were established by Qi and his coauthors in recent years.

\subsubsection{}
In~\cite[Theorem~1]{gamma-psi-batir.tex-jcam} and~\cite[Theorem~1]{gamma-psi-batir.tex-rgmia}, by virtue of the method used in~\cite[Theorem~2.4]{gamma-fun-ineq-batir} and the inequality~\eqref{log-mean-ineq} for $i=0$, the inequality~\eqref{Elezovic-Pecaric-ineq-lower} and the right-hand side inequality in~\eqref{kershaw-batir} were recovered.

\subsubsection{}
In~\cite[Theorem~2]{gamma-psi-batir.tex-jcam} and~\cite[Theorem~2]{gamma-psi-batir.tex-rgmia}, the decreasing monotonicity of the function~\eqref{bustol-ismail-AMM} and the right-hand side inequality in~\eqref{second-egp-thm4} were extended and generalized to the logarithmically complete monotonicity, and the inequality~\eqref{Elezovic-Pecaric-ineq-lower} was generalized to a decreasing monotonicity.

\begin{thm}[{\cite[Theorem~2]{gamma-psi-batir.tex-jcam} and~\cite[Theorem~2]{gamma-psi-batir.tex-rgmia}}]\label{log-complete-fcn}
For $s,t\in\mathbb{R}$ with $s\ne t$, the function
\begin{equation}\label{f_s,t}
\biggl[\frac{\Gamma(x+s)}{\Gamma(x+t)}\biggr]^{1/(s-t)}
\frac1{e^{\psi(L(s,t;x))}}
\end{equation}
is decreasing and
\begin{equation}\label{g_s,t}
\biggl[\frac{\Gamma(x+s)}{\Gamma(x+t)}\biggr]^{1/(t-s)} e^{\psi(A(s,t;x))}
\end{equation}
is logarithmically completely monotonic on $(-\min\{s,t\},\infty)$, where
$$
L(s,t;x)=L(x+s,x+t)\quad \text{and}\quad A(s,t;x)=A(x+s,x+t).
$$
\end{thm}

\subsubsection{}
In \cite{new-upper-kershaw.tex, new-upper-kershaw-JCAM.tex}, the upper bounds in~\eqref{gki2}, \eqref{second-egp-thm4}, \eqref{kershaw-batir}, \eqref{batir-psi-ineq} and related inequalities in~\cite{ratio-gamma-polynomial.tex-jcam, ratio-gamma-polynomial.tex-rgmia, gamma-batir.tex-jcam, gamma-batir.tex-rgmia} were refined and extended as follows.

\begin{thm}[\cite{new-upper-kershaw.tex, new-upper-kershaw-JCAM.tex}]\label{identric-kershaw-thm}
The inequalities
\begin{equation}\label{identric-kershaw-ineq-equiv}
\biggl[\frac{\Gamma(a)}{\Gamma(b)}\biggr]^{1/(a-b)}\le e^{\psi(I(a,b))}
\end{equation}
and
\begin{equation}\label{batir-psi-ineq-ref-equiv}
\frac{(-1)^{n}\bigl[\psi^{(n-1)}(a) -\psi^{(n-1)}(b)\bigr]}{a-b} \le(-1)^n\psi^{(n)}(I(a,b))
\end{equation}
for $a>0$ and $b>0$ hold true.
\end{thm}

\begin{rem}
The basic tools to prove~\eqref{identric-kershaw-ineq-equiv} and~\eqref{batir-psi-ineq-ref-equiv} are an inequality in~\cite{cargo} and and a complete monotonicity in \cite{subadditive-qi-3.tex} respectively. They may be recited as follows:
\begin{enumerate}
\item
If $g$ is strictly monotonic, $f$ is strictly increasing, and $f\circ g^{-1}$ is convex (or concave, respectively) on an interval $I$, then
\begin{equation}\label{carton-ineq}
g^{-1}\left(\frac1{t-s}\int_s^tg(u)\td u\right) \le
f^{-1}\left(\frac1{t-s}\int_s^tf(u)\td u\right)
\end{equation}
holds (or reverses, respectively) for $s,t\in I$. See also \cite[p.~274, Lemma~2]{bullenmean} and \cite[p.~190, Theorem~A]{f-mean}.
\item
The function
\begin{equation}\label{psi-abs-minus-cm-1}
x\bigl|\psi^{(i+1)}(x)\bigr|-\alpha\bigl|\psi^{(i)}(x)\bigr|,\quad i\in\mathbb{N}
\end{equation}
is completely monotonic on $(0,\infty)$ if and only if $0\le\alpha\le i$. See also~\cite{subadditive-qi-guo.tex}.
\end{enumerate}
\end{rem}

\begin{rem}

By the so-called G-A convex approach, the inequality~\eqref{identric-kershaw-ineq-equiv} was recovered in~\cite{Zhang-Morden}: For $b>a>0$,
\begin{equation}
[b-L(a,b)]\psi(b)+[L(a,b)-a]\psi(a)<\ln\frac{\Gamma(b)}{\Gamma(a)}<(b-a)\psi(I(a,b)).
\end{equation}
See also \href{http://www.ams.org/mathscinet-getitem?mr=2413632}{MR2413632}. Moreover, by the so-called geometrically convex method, the following double inequality was shown in~\cite[Theorem~1.2]{Xiao-Ming-JIPAM}: For positive numbers $x$ and $y$,
\begin{equation}
\frac{x^x}{y^y}\biggl(\frac{x}y\biggr)^{y[\psi(y)-\ln y]}e^{y-x}\le \frac{\Gamma(x)}{\Gamma(y)} \le\frac{x^x}{y^y}\biggl(\frac{x}y\biggr)^{x[\psi(x)-\ln x]}e^{y-x}.
\end{equation}
\end{rem}

\subsubsection{}
In~\cite{subadditive-qi-guo.tex, subadditive-qi-3.tex}, the function
\begin{equation}\label{psi-abs-minus-cm-2}
\alpha\bigl|\psi^{(i)}(x)\bigr|-x\bigl|\psi^{(i+1)}(x)\bigr|
\end{equation}
was proved to be completely monotonic on $(0,\infty)$ if and only if $\alpha\ge i+1$. Utilizing the inequality~\eqref{carton-ineq} and the completely monotonic properties of the functions~\eqref{psi-abs-minus-cm-1} and~\eqref{psi-abs-minus-cm-2} yields the following double inequality.

\begin{thm}[{\cite[Theorem~1]{new-upper-kershaw-2.tex} and~\cite[Theorem~1]{new-upper-kershaw-2.tex-mia}}] \label{new-upper-2-thm-1}
For real numbers $s>0$ and $t>0$ with $s\ne t$ and an integer $i\ge0$, the inequality
\begin{equation}\label{new-upper-main-ineq}
(-1)^i\psi^{(i)}(L_p(s,t))\le \frac{(-1)^i}{t-s}\int_s^t\psi^{(i)}(u)\td u \le(-1)^i\psi^{(i)}(L_q(s,t))
\end{equation}
holds if $p\le-i-1$ and $q\ge-i$.
\end{thm}

\begin{rem}
The double inequality~\eqref{new-upper-main-ineq} recovers, extends and refines inequalities~\eqref{Elezovic-Pecaric-ineq-lower}, \eqref{batir-psi-ineq}, \eqref{identric-kershaw-ineq-equiv} and~\eqref{batir-psi-ineq-ref-equiv}.
\end{rem}

\begin{rem}
A natural question is whether the above sufficient conditions $p\le-i-1$ and $q\ge-i$ are also necessary for the inequality~\eqref{new-upper-main-ineq} to be valid.
\end{rem}

\subsubsection{}
As generalizations of the inequalities~\eqref{Elezovic-Pecaric-ineq-lower}, \eqref{batir-psi-ineq}, the decreasing monotonicity of the function~\eqref{f_s,t}, and the left-hand side inequality in~\eqref{new-upper-main-ineq}, the following monotonic properties were presented.

\begin{thm}[{\cite[Theorem~3]{new-upper-kershaw-2.tex} and~\cite[Theorem~3]{new-upper-kershaw-2.tex-mia}}] \label{new-upper-2-thm-3}
If $i\ge0$ is an integer, $s,t\in\mathbb{R}$ with $s\ne t$, and $x>-\min\{s,t\}$, then the function
\begin{equation}\label{psi-minus-mon}
(-1)^i\left[\psi^{(i)}(L_p(s,t;x)) -\frac{1}{t-s}\int_s^t\psi^{(i)}(x+u)\td u\right]
\end{equation}
is increasing with respect to $x$ for either $p\le-(i+2)$ or $p=-(i+1)$ and decreasing with respect to $x$ for $p\ge1$, where $L_p(s,t;x)=L_p(x+s,x+t)$.
\end{thm}

\begin{rem}
It is not difficult to see that the ideal monotonic results of the function~\eqref{psi-minus-mon} should be as follows.
\end{rem}

\begin{conj}\label{new-upper-2-thm-3-conj}
Let $i\ge0$ be an integer, $s,t\in\mathbb{R}$ with $s\ne t$, and $x>-\min\{s,t\}$. Then the function~\eqref{psi-minus-mon} is increasing with respect to $x$ if and only if $p\le-(i+1)$ and decreasing with respect to $x$ if and only if $p\ge-i$.
\end{conj}

\begin{rem}
Corresponding to Conjecture~\ref{new-upper-2-thm-3-conj}, the complete monotonicity of the function~\eqref{psi-minus-mon} and its negative may also be discussed.
\end{rem}

\subsection*{Acknowledgements}
This article was ever reported on 16 February 2009 as a talk in the seminar held at the RGMIA, School of Engineering and Science, Victoria University, Australia, while the author was visiting the RGMIA between March 2008 and February 2009 by the grant from the China Scholarship Council. The author expresses thanks to Professors Pietro Cerone and Server S.~Dragomir and other local colleagues at Victoria University for their invitation and hospitality throughout this period.

\end{document}